\documentclass[10pt]{amsart}
\usepackage{amssymb}
\usepackage{MnSymbol}

\usepackage{amscd}
\usepackage{amsthm}
\usepackage[dvips]{graphicx}
\usepackage{color}
\usepackage[all]{xy}
\usepackage{verbatim}

\newtheorem{Lemma1}{{Lemma}}[section]
\newtheorem{Theo1}[Lemma1]{{Theorem}}
\newtheorem*{Theo2}{{Theorem}}
\newtheorem{Def1}[Lemma1]{{Definition}}
\newtheorem{Prop1}[Lemma1]{{Proposition}}
\newtheorem{Claim1}[Lemma1]{{Claim}}
\newtheorem{Rem1}[Lemma1]{{Remark}}
\newtheorem{Cor1}[Lemma1]{{Corollary}}
\newtheorem{Ex1}[Lemma1]{{Example}}
\newtheorem{Not1}[Lemma1]{{Notation}}

\newenvironment{Lemma}{\begin{Lemma1}}{\end{Lemma1}}
\newenvironment{Def}{\begin{Def1}\rm}{\end{Def1}}
\newenvironment{Prop}{\begin{Prop1}}{\end{Prop1}}
\newenvironment{Rem}{\begin{Rem1}\rm}{\end{Rem1}}
\newenvironment{Theorem}{\begin{Theo1}}{\end{Theo1}}

\newenvironment{Cor}{\begin{Cor1}}{\end{Cor1}}

\newenvironment{Example}{\begin{Ex1}\rm}{\end{Ex1}}

\title[DG-Hopkins-Levitzki Theorem]{DG-Semiprimary DG-Algebras, Acyclicity and 
Hopkins-Levitzki Theorem for DG-Algebras}

\author{Alexander Zimmermann}
\address{\newline
Universit\'e de Picardie,
\newline D\'epartement de Math\'ematiques et LAMFA (UMR 7352 du CNRS),
\newline 33 rue St Leu,
\newline F-80039 Amiens Cedex 1,
\newline France}
\email{alexander.zimmermann@u-picardie.fr}

\date{February 21, 2025}

\newcommand{\lra}{\longrightarrow}

\newcommand{\ra}{\rightarrow}
\newcommand{\sdp}{\times\kern-.2em\vrule height1.1ex depth-.05ex}
\newcommand{\epi}{\lra \kern-.8em\ra}

\newcommand{\N}{{\mathbb N}}

\newcommand{\ol}{\overline}

\newcommand{\Z}{{\mathbb Z}}

\newcommand{\im}{\textup{im}}
\newcommand{\dickebox}{{\vrule height5pt width5pt depth0pt}}

\newcommand{\ann}{\textup{ann}}

\newcommand{\dgrad}{\textup{dgrad}}

 \normalbaselines

\begin{document}

\begin{abstract}
We study the analogue of the Hopkins-Levitzky Theorem for dg-algebras $(A,d)$. 
We first consider the Hopkins approach. Here we show that 
for acyclic dg-algebras with graded-Artinian algebras of cycles $\ker(d)$, we also 
have that $(A,d)$ is left dg-Noetherian, and 
we show that acyclic dg-Artinian dg-algebras are dg-Noetherian. Then, 
studying the Levitzki approach,
we consider a definition of a dg-semiprimary algebra.
For dg-semiprimary dg-artinian 
dg-algebras $(A,d)$, we show that all dg-simple dg-modules are acyclic, and so are all
dg-modules with finite dg-composition length. We finally show that 
dg-Artinian dg-semiprimary  dg-algebras with nilpotent dg-radical $\dgrad_2(A,d)$  are  
dg-Noetherian and acyclic. 
\end{abstract}

\maketitle

\section*{Introduction}

A very classical result in non commutative ring theory is the Hopkins-Levitzki theorem. 
Hopkins proved in \cite{Hopkins} that a left Artinian ring is left  Noetherian. Independently and
basically the same year Levitzki showed in \cite{Levitzki} that left Artianian rings are 
semiprimary, i.e. they have a nilpotent Jacobson radical with semisimple radical quotient. He further
used this result to show that a left Artinian algebra is also left Noetherian.

A great number of research studied the ring theory of graded rings. 
For various aspects we cite   
Beattie, Chen, D\u asc\u alescu, El Fadil,  Fang, Goodearl, 
Huang, Ion, Liang, C. N\u ast\u asescu, L. N\u ast\u asescu,  
 Rio-Montes, Stafford, van den Bergh, van Geel, van Oystaeyen,
\cite{Beattie, Chen, DascalescuIon, Dascalescu,  ElFadilQuebec, GoodearlStafford, 
NastacescuvdBerghvOystaen, vanOystaenGeel}. Much of this work is presented and extended in the two 
monographs  \cite{gradedrings} and \cite{NastasescuVanOystaen}. 
In particular, in \cite[2.9.7 Corollary]{NastasescuVanOystaen} a graded version 
of Hopkins' theorem is proved. 

Differential graded algebras (dg-algebras for short) 
were defined by Cartan  \cite{Cartandg}. They proved to be a highly 
useful structure for considerations in homological algebra, in algebraic geometry, in topology, and many
other subjects. More precisely, let $R$ be a commutative ring. 
A differential graded $R$-algebra $(A,d)$ is a $\Z$-graded algebra with a degree $1$ homogeneous  
$R$-linear endomorphism $d$ of $A$ such that $d^2=0$ and satisfying 
$$d(a\cdot b)=d(a)\cdot b+(-1)^{|a|}a\cdot d(b)$$ 
for all homogeneous $a,b\in A$ where $|a|$ denotes the degree of $a$. For further undefined notations, see Section~\ref{revisitedsect} below.

Only quite recently the ring theory of this class or rings was studied. 
The first paper in this direction was Aldrich and Garcia Rozas~\cite{Tempest-Garcia-Rochas}
where acyclic algebras were classified and studied in great detail. Then Orlov \cite{Orlov1} 
studied aspects of dg-algebras which are finite dimensional over some field $R$
mostly with respect to considerations in algebraic geometry. Goodbody \cite{Goodbody} used Orlov's 
approach for a version of a dg-Nakayama's lemma. 
In \cite{dgorders} many ring theoretic general results are proved in the dg-version, 
and a version of classical orders over Dedekind domains, with additional dg-structure, as well 
as dg-projective class grops and dg-id\`eles are studied. In particular, a more suitable 
version of a dg-Nakayama lemma is shown. In \cite{dgBrauer} a dg-Brauer 
group is defined and compared with the classical Brauer group. In \cite{dgGoldie} 
Ore localisation for dg-algebras was defined 
and an analogue of Goldie's theorem in the dg-situation is proved. In \cite{dgfields} 
a dg-version of division rings were defined and completely classified. In \cite{dgseparable}
we defined a version of dg-separable dg-extensions of algebras and classified the 
(graded) commutative case completely. As one of the striking differences to the classical situation we 
gave in \cite{dgorders} an example of dg-algebra, which is simple in the sense 
that it does not allow a non trivial twosided dg-ideal, but whose category of 
dg-modules is not semisimple. In the sequel, 
the latter property will be called semisimple in the categorical sense. 
Hence, unlike the classical situation a dg-simple dg-algebra is not 
dg-semisimple in the categorical sense as is shown by Example~\ref{twobuytwomatrices}. 

In the present paper we continue this research in an attempt to see what happens with 
the Hopkins-Levitzki theorem in the dg-setting. The classical proof of the 
Hopkins-Levitzki theorem uses that an Artinian primitive algebra is simple, and simple 
algebras are semisimple in the categorical sense, where 
the latter property which is false in the dg-setting.

Nevertheless, 
the results in \cite{dgorders, dgGoldie, dgfields, dgseparable} suggest that 
graded conditions on the graded algebra $\ker(d)$ of cycles are likely to imply 
the conclusions for the dg-algebra $(A,d)$. 
So, we prove that if $\ker(d)$ is graded-Noetherian (resp. graded-Artinian) 
then $(A,d)$ is dg-Noetherian (resp. dg-Artinian). As a consequence, using the 
graded Hopkins-Levitzki theorem we get that if $\ker(d)$ is graded-Artinian then $(A,d)$ is 
dg-Artinian and dg-Noetherian (cf Corollary~\ref{gr-artinianimpliesdgartanddgnoeth} below). 
Further, for acyclic dg-algebras we simply get the 
plain Hopkins-Levitzki result, namely dg-Artinian acyclic dg-algebras are 
dg-Noetherian (cf Theorem \ref{acyclicHopkins} below).
We then pass to the Levitzki approach and define dg-semiprimary dg-algebras in a generalized 
version motivated by \cite{CR}, replacing the nilpotence of the Jacobson radical by a lifting
property of idempotents. We then show that this gives algebras which are acyclic 
in case they have finite composition length (cf Theorem~\ref{alldgmodulesareacyclic} 
in connection with \ref{dgsemiprimaryhomologyofmaximals} below). 
We finally prove in Theorem~\ref{dgsemiprimaryLevitzki}  that left dg-Artinian left dg-semiprimary 
dg-algebras with nilpotent $\dgrad_2(A)$ are left dg-Noetherian and acyclic.  

The paper is organised as follows. In Section~\ref{revisitedsect} we recall the 
necessary results and definitions needed in this paper. Section~\ref{gradedlevitzkiconsequences}
draws consequences from the graded version of the Hopkins-Levitzki 
Theorem~\ref{gradedHopkins}. In particular, it is shown
in Theorem~\ref{acyclicHopkins} the dg-version in case of acyclic algebras.
Section~\ref{dgprimarysection} defines dg-semiprimary algebras and studies their properties,
in particular with respect to their homology. Theorem~\ref{alldgmodulesareacyclic} in connection
with Proposition~\ref{dgsemiprimaryhomologyofmaximals} provides a quite definite picture. 
Section~\ref{thefinalresult} then shows the final result for a Hopkins-Levitzki Theorem in the 
dg-semiprimary dg-algebra setting.

\subsection*{Acknowledgement:}
I wish to thank Raphael Bennett-Tennenhaus for very instructive discussions 
in May 2024 on the subject.

\section{Dg-algebras and their dg-modules revisited}

\label{revisitedsect}

In this section we shall recall definitions from the general theory, 
following mainly \cite{Yekutielibook}, and results proved earlier in \cite{dgorders}, as
far as they are needed in this paper.  

\subsection*{dg-algebras and dg-modules:}
Let $R$ be a commutative ring. 
A differential graded $R$-algebra $(A,d)$ is a $\Z$-graded 
$R$-algebra $A$ together with an $R$-linear endomorphism $d$ of $A$
which is homogeneous of degree $1$, such that $d^2=0$ and 
$$d(a\cdot b)=d(a)\cdot b+(-1)^{|a|}a\cdot d(b)$$ 
for all homogeneous $a,b\in A$, where $|a|$ denotes the degree of $a$.
A differential graded left module $(M,\delta)$ over the differential graded 
$R$-algebra $(A,d)$ is a $\Z$-graded left $A$-module $M$
and an $R$-linear endomorphism $\delta$ of $M$ of degree $1$
with $\delta^2=0$ satisfying 
$$ 
\delta(a\cdot m)=d(a)\cdot m+(-1)^{|a|}a\cdot \delta(m)
$$
for all homogeneous $a\in A$ and $m\in M$. 
A dg-right module is  a $\Z$-graded right $A$-module $M$
and an $R$-linear endomorphism $\delta$ of $M$ of degree $1$
with $\delta^2=0$ satisfying 
$$ 
\delta(m\cdot a)=\delta(m)\cdot a+(-1)^{|m|}m\cdot d(a)
$$
for all homogeneous $a\in A$ and $m\in M$. 
We abbreviate `differential graded' by `dg' for short. 
If the differential $d$ (resp. $\delta$) is $0$, we use 
the notion `graded-' instead of `differential graded'.

\subsection*{Homomorphisms between dg-algebras and between dg-modules:}
Given two dg-modules $(M,\delta_M)$ and $(N,\delta_N)$ over a dg-algebra 
$(A,d)$ we say that a homomorphism of dg-modules 
$(M,\delta_M)\ra(N,\delta_N)$ is a degree $0$ homogeneous 
$A$-module homomorphism $\varphi:M\rightarrow N$ such that 
$\varphi\circ\delta_M=\delta_N\circ\varphi$. Similarly, for two
dg-algebras $(A,d_A)$ and $(B,d_B)$ we say that a homomorphism of dg-algebras
$\alpha:(A,d_A)\ra (B,d_B)$ is a degree $0$ homogeneous graded algebra
homomorphism $A\ra B$ such that $\alpha\circ d_A=d_B\circ\alpha$. 
Isomorphisms, monomorphisms and epimorphisms are defined in the 
corresponding category of dg-modules resp. dg-algebras. 

\subsection*{dg-Noetherian and dg-Artinian dg-modules:}
For a dg-module $(M,\delta)$ over $(A,d)$ a dg-sub\-module $(N,\delta)$ of $(M,\delta)$ is 
a graded submodule $N$ of $M$ such that $\delta(N)\subseteq N$. 
In particular, $(N,\delta|_N)$ is a dg-module over $(A,d)$. 
A dg-module $(M,\delta)$ is dg-Noetherian (resp. dg-Artinian)
if it satisfies the ascending (resp. descending) chain condition on 
dg-submodules. A dg-module $(S,\delta)$ is dg-simple if it does not allow a 
dg-submodule other than $0$ or $S$. A dg-algebra $(A,d)$ is dg-simple if 
it is dg-simple as dg-bimodule over $(A,d)$. 

\subsection*{dg-radical and dg-Nakayamas lemma:}
Recall that we introduced the dg-radical of a dg-algebra in \cite{dgorders}
as 
$$\dgrad_2(A)=\bigcap_{(S,\delta)\textup{ dg-simple dg-module}}\ann(S,\delta),$$
where $\ann$ is the left annihilator in $A$, and it was proved in \cite{dgorders} 
\begin{itemize}
\item that if $(M,\delta)$ is a dg-module, finitely generated as a 
dg-module, then $$\left(\dgrad_2(A)\cdot M=M\right)\Rightarrow M=0.$$ 
\item 
Further, it was proved in \cite{dgorders} that $(M,\delta)$ is 
dg-Noetherian if and only if all 
dg-submodules of $(M,\delta)$ are finitely generated.
\end{itemize} 

As a consequence, if $(A,d)$ is dg-Noetherian as a left dg-module over 
itself, then for all integers $k>0$ the left dg-ideal
$(\dgrad_2(A)^k$ is finitely generated as a left dg-module. 
Hence, the sequence of twosided dg-ideals $(\dgrad_2(A))^k$ for $k\in\N$ is 
decreasing, and if $(A,d)$ is dg-Noetherian as a 
dg-left module over itself, also 
finitely generated. Hence the sequence has to be stationary at some point, and 
by the dg-Nakayama Lemma, we get that $\dgrad_2(A)$ is nilpotent. 

\subsection*{dg-semisimplicity:}
We need to recall that we have two concepts of semisimplicity. 

The first concept is what we call the categorical semisimplicity. 
A dg-algebra $(A,d)$ is called semisimple in the categorical sense if 
the category of dg-modules is a semisimple category, i.e. 
and short exact sequence of dg-modules is split. 

The second concept is semisimplicity on the Jacobson sense. Here we say that a
dg-module $(M,\delta)$ is semisimple in the Jacobson sense if 
$(M,\delta)$ is a direct sum of dg-simple dg-submodules 
$(S_i,\delta)\leq (M,\delta)$ (for $i\in I$, an index set):
$$(M,\delta)=\bigoplus_{i\in I}(S_i,\delta).$$
An algebra is called dg-semisimple in the Jacobson sense if the 
regular dg-module is dg-semisimple in the Jacobson sense. 

Trivially, a semisimple algebra in the categorical sense is 
semisimple in the Jacobson sense as well. 
For an ungraded Artinian algebra, categorical semisimplicity 
and Jacobson semisimplicity coincide. However, for dg-algebras 
this not true. Examples were given in \cite[Example 4.35]{dgorders},
which re-appears as Example~\ref{twobuytwomatrices} below.

\section{Consequences of the graded Hopkins-Levitzki Theorem}

\label{gradedlevitzkiconsequences}

Recall the graded Hopkins-Levitzki Theorem.

\begin{Theorem} \cite[2.9.7 Corollary]{NastasescuVanOystaen} \label{gradedHopkins}
Let $A$ be a graded algebra. If $A$ is graded-Artinian, then $A$ is graded-Noetherian. 
\end{Theorem}

\begin{Prop}\label{kerdimpliesA}
Let $(A,d)$ be a dg-algebra and let $I\leq J$. If 
$I\cap\ker(d)=J\cap\ker(d)$ and $d(I)=d(J)$, then 
$I=J$. 
\end{Prop}

Proof. 
Since $I\leq J$, we get $d(I)\leq d(J)$ and $I\cap\ker(d)\leq J\cap\ker(d)$. 
Hence, the following diagram with exact rows and columns is commutative.
$$
\xymatrix{&0\ar[d]&0\ar[d]&0\ar[d]\\
0\ar[r]&I\cap\ker(d)\ar[r]\ar[d]&I\ar[r]^d\ar[d]&d(I)\ar[r]\ar[d]&0\\
0\ar[r]&J\cap\ker(d)\ar[r]\ar[d]&J\ar[r]^d\ar[d]&d(J)\ar[r]\ar[d]&0\\
&0&J/I\ar[d]&0\\
&&0
}
$$
By the snake lemma we also get that the lowest row is commutative:
$$
\xymatrix{&0\ar[d]&0\ar[d]&0\ar[d]\\
0\ar[r]&I\cap\ker(d)\ar[r]\ar[d]&I\ar[r]^d\ar[d]&d(I)\ar[r]\ar[d]&0\\
0\ar[r]&J\cap\ker(d)\ar[r]\ar[d]&J\ar[r]^d\ar[d]&d(J)\ar[r]\ar[d]&0\\
&0\ar@{-->}[r]&J/I\ar[d]\ar@{-->}[r]&0\\
&&0
}
$$
This shows the statement. \dickebox

\begin{Cor} \label{grnoethimpliesdgnoeth}
Let $(A,d)$ be a dg-algebra. 
\begin{itemize}
\item If $\ker(d)$ is graded-Noetherian, then $(A,d)$ 
is dg-Noetherian. 
\item If $\ker(d)$ is graded-Artinian, then $(A,d)$ 
is dg-Artinian.
\end{itemize}
\end{Cor}

Proof. 
First note that if $I$ is a dg-ideal of $A$, then $I\cap\ker(d)$ is a graded ideal of $\ker(d)$.
Indeed, if $y\in\ker(d)$ and $x\in A$, then $d(x)\in\ker(d)$, but also $d(x)\in I$, and $I$ is an ideal. Hence 
$y\cdot d(x)\in I$, which shows
$y\cdot d(x)\in\ker(d)\cap I$.  

\medskip

Further, for a dg-ideal $I$ of $A$ we get that $d(I)$ is an ideal of $\ker(d)$. 
Indeed, clearly $I$ is closed under addition. Let $y\in\ker(d)$ be homogeneous and $x\in I$. 
Then $$y\cdot d(x)=(-1)^{|y|}d(yx)-(-1)^{|y|}d(y)\cdot x=(-1)^{|y|}d(yx)\in d(I).$$
Hence $d(I)\leq\ker(d)$. 

\medskip

Let $$I_1\leq I_2\leq\dots\leq A$$ be an increasing sequence of dg-ideals of $A$, then 
$$I_1\cap\ker(d)\leq I_2\cap\ker(d)\leq\dots\leq\ker(d)$$
and 
$$d(I_1)\leq d(I_2)\leq \dots\leq \ker(d)$$
are increasing sequences of graded ideals of $\ker(d)$. 
If $\ker(d)$ is graded-Noetherian, then there is $n$ with $I_n\cap\ker(d)=I_{n+k}\cap\ker(d)$
for all $k\in\N$ and there is $m$ with $d(I_m)=d(I_{m+k})$  
for all $k\in\N$. For $s:=max(n,m)$ we hence get that $d(I_s)=d(I_{s+k})$ and 
$I_s\cap\ker(d)=I_{s+k}\cap\ker(d)$ for all $k\in\N$. By  Proposition~\ref{kerdimpliesA}
we get that $I_s=I_{s+k}$ for all $k\in\N$. This shows the statement for the Noetherian property. 

The statement for the Artinian case follows by the very same argument, considering 
decreasing sequences of dg-ideals. \dickebox

\begin{Prop}\label{gr-artinianimpliesdgartanddgnoeth}
Let $(A,d)$ be a dg-algebra. If $\ker(d)$ is graded-Artinian, then $(A,d)$ is dg-Noetherian and dg-Artinian. 
\end{Prop}

Proof. If $\ker(d)$ is graded-Artinian, then $(A,d)$ is dg-Artinian by the second statement of 
Corollary~\ref{grnoethimpliesdgnoeth}.
If $\ker(d)$ is graded-artinina, then $\ker(d)$ is graded-Noetherian by Theorem~\ref{gradedHopkins}. 
Then, using the first statement of Corollary~\ref{grnoethimpliesdgnoeth}, we see that $(A,d)$ is dg-Noetherian. \dickebox

\begin{Rem} \label{diagramgeneral}
To summarize we have shown what is displayed in the following diagram:
$$\xymatrix{
\ker(d)\textup{ gr-Artinian}\ar@{=>}[r]^{\textup{Cor~\ref{grnoethimpliesdgnoeth}}}\ar@{=>}[d]_{\textup{Thm~\ref{gradedHopkins}}}
&(A,d)\textup{ dg-Artinian}\\
 \ker(d)\textup{ gr-Noetherian}\ar@{=>}[r]^{\textup{Cor~\ref{grnoethimpliesdgnoeth}}}&(A,d)\textup{ dg-Noetherian}
}$$
\end{Rem}

\medskip

In order to establish a dg-Hopkins' Theorem we need to show that if $(A,d)$ is dg-Artinian, then $\ker(d)$ is graded-Artinian as well. 

\begin{Lemma}\label{fgidealsofcycles}
Let $(A,d)$ be an acyclic dg-algebra. Let $\ol I$ and $\ol J$ be graded
ideals of $\ker(d)$ such that $\ol I\leq\ol J$.  Then 
$$\ol I\neq\ol J\Rightarrow A\cdot\ol I\neq A\cdot\ol J.$$
\end{Lemma}

Proof. 
Since the ideals $\ol I$ and $\ol J$ are 
graded  $\ker(d)$-ideals, we may express
$\ol I=\sum_i\ker(d)\cdot x_i$ and $\ol J=\sum_j\ker(d)y_j$
for elements $x_i,y_j\in\ker(d)$. 
Then $I:=A\cdot\ol I=\sum_iAx_i$ and $J:=A\cdot \ol J=\sum_jAy_j$.
Hence, using that $A$ is acyclic, and therefore $\ker(d)=d(A)$, 
we get $d(I)=\sum_id(A)x_i=\sum_i\ker(d)x_i=\ol I$. 
Similarly, $d(J)= \ol J$. 
Hence, $I=J$ implies $\ol I= d(I)=d(J)=\ol J$. \dickebox

\begin{Theorem}\label{acyclicHopkins}
If $(A,d)$ is a dg-Artinian acyclic dg-algebra, then $(A,d)$ is dg-Noetherian. 
\end{Theorem}

Proof.
Indeed, By Lemma~\ref{fgidealsofcycles} we see that $\ker(d)$ is graded-Artinian. By 
Theorem~\ref{gradedHopkins} we get that $\ker(d)$ is graded-Noetherian. 
By Corollary~\ref{grnoethimpliesdgnoeth}  we see that $(A,d)$ is dg-Noetherian. \dickebox

\begin{Rem}
Lemma~\ref{fgidealsofcycles} follows also from \cite{Tempest-Garcia-Rochas} 
since there it is shown that for acyclic dg-algebras
$A\otimes_{\ker(d)}-$ is an equivalence between the category of graded modules over $\ker(d)$ and 
dg-modules over $(A,d)$.  We believe however that our short proof is simpler. 
\end{Rem}

\begin{Rem}
{\em If $(A,d)$ is acyclic,} then the diagram from Remark~\ref{diagramgeneral} becomes
$$\xymatrix{
\ker(d)\textup{ gr-Artinian}\ar@{=>}@/^/[r]^{\textup{Cor~\ref{grnoethimpliesdgnoeth}}}
\ar@{=>}[d]_{\textup{Thm~\ref{gradedHopkins}}}
&(A,d)\textup{ dg-Artinian}\ar@{=>}@/^/[l]^{\textup{Lemma~\ref{fgidealsofcycles}}}\\
 \ker(d)\textup{ gr-Noetherian}\ar@{=>}[r]^{\textup{Cor~\ref{grnoethimpliesdgnoeth}}}&(A,d)\textup{ dg-Noetherian}
}$$
\end{Rem}

\begin{Rem} A possible counterexample for the general case would be a dg-algebra $(A,d)$, 
where $\ker(d)$ is not graded-Artinian, but $(A,d)$ is dg-Artinian.
\end{Rem}

\section{Dg-semiprimary dg-algebras}

\label{dgprimarysection}

The classical Levitzki-theorem uses semiprimary algebras. We shall give a dg-version
of this concept.

\begin{Def}
Let $(A,d)$ be a dg-algebra. We say that $(A,d)$ is {\em left dg-semiprimary} if 
 $A/\dgrad_2(A)$ is left 
dg-semisimple in the categorical sense, i.e. the category of left
$A/\dgrad_2(A)$-modules is a semisimple category, 
and if any homogeneous idempotent in $A/\dgrad_2(A)$ 
is the image of a homogeneous idempotent in
$A$ under the natural map. 
\end{Def}

\begin{Rem} \label{huchitsacyclic}
\begin{itemize}
\item
Note that the condition on lifting idempotents holds 
in particular if $\dgrad_2(A)$ is nilpotent. Mostly, one asks this latter 
condition, but here we shall follow \cite{CR} where an analogue of this definition is 
proposed.
\item
Note further that a homogeneous idempotent is necessarily in degree $0$. The 
idempotent lifting property statement is hence a question on the degree $0$ 
component of $A$. 
\item
We see that if $(A,d)$ left dg-semiprimary, then 
$A/\dgrad_2(A)$ is a direct sum of left dg-simple dg-modules.
By \cite[Theorem 5.3]{Tempest-Garcia-Rochas} we further get that
this is equivalent with $(A/\dgrad_2(A),\ol d)$ is acyclic,
which in turn is equivalent with $1\in \im(\ol d)$. This then implies that
$$A/\dgrad_2(A)=\ker(\ol d)\oplus \ker(\ol d)\cdot\ol y=
\ker(\ol d)\oplus \ol y\cdot\ker(\ol d)$$ 
for some element $\ol y$ with $\ol d(\ol y)=1$. 
\end{itemize}
\end{Rem}

\begin{Example} \label{twobuytwomatrices}
For a field $K$ let $A=Mat_2(K)$ with differential being 
$$d(\left(\begin{array}{cc}x&y\\z&w\end{array}\right))=
\left(\begin{array}{cc}z&w-x\\0&z\end{array}\right)$$
is a dg-algebra with a unique dg-simple left dg-module 
(cf \cite[Example 3.6 and Example 4.35]{dgorders}). 
Its left annihilator is $0$, and hence $\dgrad_2(A)=0$.
Nevertheless, $A$ is simple Artinian, whence dg-simple dg-Artinian, but not 
dg-semisimple in the categorical sense, and not even in the Jacobson sense. 
Therefore, $(A,d)$ is not 
left dg-semiprimary. However, $(A,d)$ is acyclic, but $\ker(d)$ is not graded-semisimple.
\end{Example}

\begin{Example}
$(A,d)=(K[X]/X^2,d)$ with $d(X)=1$ is an acyclic dg-algebra. Also $\dgrad_2(A)=0$
and $\ker(d)=K$ is graded semisimple. Hence $(A,d)$ is dg-semiprimary.  
Actually all dg-division algebras \cite{dgfields} are dg-semiprimary.
\end{Example}

\begin{Rem}
If $(A,d)$ is dg-semiprimary left dg-Artinian, then the left dg-module  
$A/\dgrad_2(A)$ over the dg-algebra $(A,d)$ is a 
finite direct sum of dg-simple left dg-modules over $(A,d)$.

However, $(A,d)$ semiprimary implies that $(A/\dgrad_2(A),\ol d)$ 
is a dg-semisimple dg-module over  $(A/\dgrad_2(A),\ol d)$ in the categorical sense,  which
is a stronger property. 
\end{Rem}

\begin{Example}
Let $A$ be a finite dimensional semisimple $K$-algebra over a field $K$. This is a dg-algebra,
considering the trivial grading on $A$ and the differential $0$. Then $\dgrad_2(A)=\textup{rad}(A)=0$
and the category of $A$-modules is semisimple. However, the homology of $(A,0)$ is simply $A$, 
which is not $0$. This is not a contradiction to what was said before in Remark~\ref{huchitsacyclic}, since 
$(A,0)-\textup{dgmod}$, the category of dg-modules over $(A,0)$ is not semisimple, as
the category of dg-modules over an algebra is the category of complexes over $A$, and not the module
category over $A$.  

Hence, as a consequence, a semiprimary algebra considered as dg-algebra concentrated in degree $0$ with 
differential $0$ is not dg-semiprimary. 
\end{Example}

\begin{Lemma} \cite[Lemma 4.19]{dgorders}\label{dgsimplesarequotientsofA}
Let $(A,d)$ be a dg-algebra and suppose that $(S,\delta)$ is a dg-simple dg-module. 
Then there is a dg-maximal dg-left ideal $({\mathfrak m},d)$ of $(A,d)$ and an integer $n\in\Z$ 
such that $(S,\delta)\simeq(A/{\mathfrak m},\ol d)[n]$.
\end{Lemma}

\begin{Lemma} \label{jacobsondgsemieimpleimpliesallsimplesareinA}
Let $(A,d)$ be a dg-algebra, and suppose that the dg-left module $(A,d)$ over $(A,d)$ 
is a finite direct sum of dg-simple dg-modules over $(A,d)$. Then for
any dg-simple dg-module $(S,\delta)$ over $(A,d)$ there is an integer $n\in\Z$ such that 
$(S,\delta)[n]$ is a direct factor of $(A,d)$. 
\end{Lemma}

Proof. 
By Lemma~\ref{dgsimplesarequotientsofA} we may suppose that 
$(S,\delta)\simeq(A/{\mathfrak m},\ol d)$ for some dg-maximal left dg-ideal $({\mathfrak m},d)$ of $(A,d)$. 
By \cite[Lemma 4.20]{dgorders} we see that there is a left dg-ideal $(T,d)$
of $(A,d)$ such that ${\mathfrak m}\oplus T=A$. But then 
$$(S,\delta)\simeq (A/{\mathfrak m},\ol d)=(({\mathfrak m}\oplus T)/{\mathfrak m},\ol d)\simeq (T,d)$$
as left dg-modules. This shows the statement. \dickebox

\begin{Lemma} \label{acyclicquotienthomologyofideal}
Let $(A,d)$ be a dg-algebra. Suppose that $(I,d)$ is a dg-ideal 
such that $(A/I,\ol d)$ is acyclic. Then the 
inclusion $I\subseteq A$ induces an isomorphism 
$H(A,d)\simeq H(I,d)$. 
\end{Lemma}

Proof. Indeed, 
we have a short exact sequence
$$
0\lra (I,d)\lra (A,d)\lra (A/I,\ol d)\lra 0
$$
of complexes. Taking homology, we get a long exact sequence
$$
\dots\lra H_{n-1}(A/I,\ol d)\lra H_n(I,d)\lra H_n(A,d)\lra H_n(A/I,\ol d)\lra \dots.
$$
Since the left and the right hand terms are $0$ by hypothesis, we get 
$H_n(A,d)\simeq H_n(I,d)$ for any $n$. 
 \dickebox

\medskip

We shall consider the consequences of Lemma~\ref{acyclicquotienthomologyofideal} for
the case of dg-semiprimary dg-algebras.

\begin{Prop} \label{dgsemiprimaryhomologyofmaximals}
Let $(A,d)$ be a left dg-semiprimary dg-algebra.
\begin{itemize}
\item 
Then any dg-simple left dg-module $(S,\delta)$ is acyclic. 
\item
Further, for any dg-maximal left dg-ideal $(M,d)$ of $(A,d)$
we get that $H(M,d)\simeq H(A,d)$.
\item
If $(\dgrad_2(A),d)$ is an acyclic dg-ideal then $(A,d)$ is acyclic.
\end{itemize}
\end{Prop}

Proof. 
Since $(A,d)$ is left dg-semiprimary, the dg-algebra 
$(A/\dgrad_2(A),\ol d)$ is semisimple in the categorical sense, 
and hence, by \cite{Tempest-Garcia-Rochas}
the algebra $(A/\dgrad_2(A),\ol d)$ is acyclic. 

Let $(S,\delta)$ be a dg-simple dg-module. 
Since 
$$\dgrad_2(A,d)=\bigcap_{(T,\partial)\textup{ dg-simple dg-module}}\ann(T,\partial),$$
we get that $\dgrad_2(A)\cdot(S,\delta)=0$, and hence 
$(S,\delta)$ is a dg-simple dg-module over the acyclic dg-algebra 
$(A/\dgrad_2(A),\ol d)$. Moreover, by 
\cite[Theorem 4.7]{Tempest-Garcia-Rochas} we get that any dg-module over an 
acyclic dg-algebra is acyclic. Hence $(S,\delta)$ is acyclic. 

Let now $(M,d)$ be a dg-maximal left dg-ideal of $(A,d)$. Then 
$(A/M,\ol d)$ is a dg-simple left dg-module over $(A,d)$. Hence 
$(A/M,\ol d)$ is acyclic by the previous statement, and 
by Lemma~\ref{acyclicquotienthomologyofideal} we get that 
the inclusion induces an isomorphism $H(M,d)\simeq H(A,d)$.

In case $(\dgrad_2(A),d)$ is an acyclic dg-ideal,  we get the statement by
applying Lemma~\ref{acyclicquotienthomologyofideal}.
\dickebox

\begin{Theorem} \label{alldgmodulesareacyclic}
Let $(A,d)$ be a dg-algebra such that every dg-simple left dg-module is acyclic. Then 
any left dg-module with finite dg-composition length is acyclic. 
In particular, if $(A,d)$ has finite composition length and all 
dg-simple left dg-modules are acyclic, then $(A,d)$ is acyclic. 
\end{Theorem}

Proof. 
We shall use induction on the composition length $n$.

For $n=1$, the module is dg-simple and we are done by hypothesis. 

Let $(M,\delta_M)$ be a left dg-module with composition length $n$. 
Hence, there is a dg-simple dg-submodule $(S,\delta_S)$ such that 
$$
0\lra (S,\delta_S)\lra (M,\delta_M)\lra (M/S,\ol\delta_M)\lra 0
$$
is a short exact sequence of left dg-modules over $(A,d)$. 
Now, the dg-composition length of $(M/S,\ol\delta_M)$ is $n-1$.
This induces a 
long exact sequence on homology part of which is
$$
\dots \lra H_k(S,\delta_S)\lra H_k(M,\delta_M)\lra H_{k}((M/S,\ol\delta_M))\lra\dots.
$$
The left term is $0$ by hypothesis, and the right term is $0$ by the induction hypothesis, using that the 
composition length of $(M/S,\ol\delta_M)$ is $n-1$. Hence $H_k(M,\delta_M)=0$. 
This proves the statement. \dickebox

\section{Levitzki-Hopkins Theorem for dg-semiprimary dg-algebras}

\label{thefinalresult}

We are now able to prove our final result. 

\begin{Theorem}\label{dgsemiprimaryLevitzki}
Let $(A,d)$ be a left dg-semiprimary left dg-Artinian dg-algebra
with nilpotent dg-radical $\dgrad_2(A)$. 
Then $(A,d)$ is left dg-Noetherian and acyclic.  
\end{Theorem}

\begin{Rem}
If $\dgrad_2(A)$ is finitely generated as left dg-ideal,
then for each integer $k>0$ 
also $\dgrad_2(A)^k$ is finitely generated, with generators being
words of length $k$ of these generators. Hence, by the dg-Nakayama lemma
if $$\dgrad_2(A)^n=\dgrad_2(A)^{n+1}=\dgrad_2(A)\cdot \dgrad_2(A)^n,$$
we get $\dgrad_2(A)^n=0$.  This implies that if $\dgrad_2(A)$ is finitely generated, then 
either   $\dgrad_2(A)$ is nilpotent, or the descending sequence 
of dg-ideals 
$$A\supsetneq  \dgrad_2(A)\supsetneq  \dgrad_2(A)^2 \supsetneq  \dgrad_2(A)^3 \supsetneq\dots$$
is an infinite strictly descreasing sequence of twosided dg-ideals.
\end{Rem}

Proof of Theorem~\ref{dgsemiprimaryLevitzki}. 
As $(A,d)$ is dg-Artinian dg-semiprimary with nilpotent $\dgrad_2(A)$, 
each of the dg-modules $\dgrad_2^i(A)/\dgrad_2^{i+1}(A)$ is actually a
differential graded
$A/\dgrad_2(A)$ module, and only finitely many of these quotients, 
say only for the indices $i\in\{0,\dots,n\}$, are non zero. 
$(A,d)$ dg-Artinian implies that these modules need to
be finitely generated as dg-modules, and hence $\dgrad_2^i(A)/\dgrad_2^{i+1}(A)$
is a direct sum of dg-simple dg-modules over $(A,d)$.
Therefore, for any integer $i\in\N$ the quotient 
$\dgrad_2^i(A)/\dgrad_2^{i+1}(A)$ is a direct sum of 
a finite number $k_i\in\N$  of dg-simple dg-modules over $(A,d)$ (actually over $(A/\dgrad_2(A))$).
Hence $A$ allows a dg-composition series with $\ell(A,d):=\sum_{i=0}^nk_i<\infty$ 
composition factors.
By the Zassenhaus butterfly lemma any two dg-composition series 
have the same set (counted with multiplicities) 
of dg-isomorphism classes of dg-simple dg-modules. 

Now, any increasing sequence of dg-ideals can be refined using intersections with the 
composition series into a composition series itself, and Zassenhaus' butterfly lemma 
then shows that the number of dg-composition factors equals the finite number
$\ell(A,d)$ from above. Hence, $(A,d)$ is Noetherian. 

By Proposition~\ref{dgsemiprimaryhomologyofmaximals} and 
Theorem~\ref{alldgmodulesareacyclic} we get that $(A,d)$ is acyclic. 
 \dickebox

%
%

\begin{Rem}
Theorem~\ref{dgsemiprimaryLevitzki} shows that dg-semiprimary dg-algebras with nilpotent $\dgrad_2(A)$ 
satisfying a dg-version of the Hopkins-Levitzki Theorem are necessarily acyclic. They hence fall under the 
hypotheses of Theorem~\ref{acyclicHopkins}. 
\end{Rem}

\begin{Rem}
In general an Artinian algebra is semiprimary. In order to prove this
one proceeds by first observing that quotients of left Artinian 
semiprimary algebras modulo the Jacobson radical are left Artinian left sempiprimitive.
Then one shows that Artinian left primitive algebras are simple.
Finally one uses that Artinian simple algebras are semisimple in the categorical sense.

However, Example~\ref{twobuytwomatrices} shows that  dg-simple algebras 
are not necessarily dg-semisimple in the categorical sense, and are not 
dg-semisimple in the Jacobson sense. Hence, the proof 
used in the algebra case
collapses in the differential graded case.
\end{Rem}

\end{document}